\documentclass[a4paper,fleqn]{article}
\usepackage{amsmath}
\usepackage{amssymb}
\usepackage{amsthm}
\usepackage{changebar}
\usepackage{enumerate}
\usepackage{cite,amsfonts}

\newtheorem{theorem}{Theorem}
\newtheorem{corollary}[theorem]{Corollary}

\newtheorem{proposition}[theorem]{Proposition}
\newtheorem{lemma}[theorem]{Lemma}
\theoremstyle{definition}
\newtheorem{definition}[theorem]{Definition}
\theoremstyle{remark}
\newtheorem{remark}[theorem]{Remark}
\newtheorem{acknowledgement}[theorem]{Acknowledgement}

\newcommand{\ud}{\mathrm{d}}
\newcommand{\pd}{\partial}
\newcommand{\e}{\varepsilon}
\newcommand{\bC}{\mathbb{C}}
\newcommand{\bE}{\mathbb{E}}
\newcommand{\bN}{\mathbb{N}}
\newcommand{\bR}{\mathbb{R}}
\newcommand{\cA}{\mathcal{A}}
\newcommand{\ccD}{\mathcal{D}}
\newcommand{\cE}{\mathcal{E}}
\newcommand{\cG}{\mathcal{G}}

\newcommand{\cN}{\mathcal{N}}
\newcommand{\csub}{\subset \subset}
\newcommand{\coleq}{\mathrel{\mathop:}=}
\newcommand{\coliff}{\mathrel{\mathop:}\Leftrightarrow}
\newcommand{\norm}[1]{\left\lVert#1\right\rVert}
\newcommand{\abso}[1]{\left\lvert#1\right\rvert}

\DeclareMathOperator{\supp}{supp}
\DeclareMathOperator{\dist}{dist}
\DeclareMathOperator{\pr}{pr}

\begin{document}


\title{Point value characterizations and related results in the full Colombeau algebras $\cG^e(\Omega)$ and $\cG^d(\Omega)$}

\author{Eduard Nigsch\footnote{University of Vienna, Faculty of Mathematics, Nordbergstr.~15, A-1090 Vienna, Austria, eduard.nigsch@univie.ac.at}}
\date{}
\maketitle

\begin{abstract}
We present a point value characterization for elements of the elementary full Colombeau algebra $\cG^e(\Omega)$ and the diffeomorphism invariant full Colombeau algebra $\cG^d(\Omega)$. Moreover, several results from the special algebra $\cG^s(\Omega)$ about generalized numbers and invertibility are extended to the elementary full algebra.
\end{abstract}

\vskip 1em



\section{Introduction}

Colombeau algebras \cite{ColElem} are spaces of generalized functions which serve to extend the theory of Schwartz distributions such that these can be multiplied, circumventing the well-known impossibility result by Schwartz \cite{Schwartz}. These commutative and associative differential algebras provide an embedding of the space of distributions as a linear subspace and the space of smooth functions as a faithful subalgebra.

For Schwartz distributions a concept of point values was introduced \cite{Lojasiewicz}, but an arbitrary distribution need not have a point value in this sense at every point. Furthermore, it is not possible to characterize distributions by their point values. Colombeau-type algebras of generalized functions are usually constructed as nets of smooth functions, which means that a given point can be inserted into each component of the net in order to give a generalized point value. This is not sufficient for uniquely characterizing a generalized function, though: there exist nonzero generalized functions that evaluate to zero at every classical point. However, with the introduction of generalized points one can obtain a point value characterization theorem. Note that for holomorphic generalized functions a stronger results holds, which states that such a function is zero already if its zero set has positive measure \cite{Scarpi}. Point values for Colombeau generalized functions were first introduced for $\cG^s(\Omega)$, the special Colombeau algebra on an open set $\Omega \subseteq \bR^n$ \cite{ObeChar}, and later on also for the special Colombeau algebra on a manifold \cite{foundgeom}. In the context of $p$-adic Colombeau-Egorov type generalized functions it was first claimed that classical points suffice to characterize a function \cite{Shelkovich}, but this claim was shown to be invalid later on and a characterization using generalized points was given in \cite{Eberhard}.

The aim of the present work is to introduce generalized points, numbers, and point values for the elementary full algebra $\cG^e(\Omega)$ \cite{ColElem} and the diffeomorphism invariant full algebra $\cG^d(\Omega)$ \cite{found}. Both algebras are presented in a unifying framework in \cite{GKOS}. Our main result is a point value characterization theorem for each algebra (Theorems \ref{gepvchar} and \ref{gdpvchar}) which states that two generalized functions are equal if and only if they have the same generalized point value at all generalized points.

Let us mention some applications generalized numbers and point values have found so far. First, when one does Lie group analysis of differential equations in generalized function spaces, point values allow to transfer the classical procedure for computing symmetries to the generalized case \cite{GroupAnalysis}. Second, consider mappings from the space of generalized points into the space of generalized numbers. For such mappings a discontinuous differential calculus was constructed, featuring a fundamental theorem of calculus, notions of sub-linear, holomorphic, and analytic mappings, generalized manifolds, and related results \cite{DiscontCalculus}. Using point values, elements of $\cG^s$ can be regarded as such mappings. From this viewpoint their local properties can be analyzed \cite{LocalProperties}. Moreover, point values have repeatedly turned out to be indispensable tools for doing analysis in algebras of generalized functions (cf., e.g., \cite{claudia1, claudia2, hans, scarpi2}).

\section{Preliminaries}

The number $n \in \bN$ will always denote the dimension of the underlying space $\bR^n$. $\pd \Omega$ denotes the topological boundary of a set $\Omega$. For $A \subseteq \bR^n$ we write $K \csub A$ if $K$ is a compact subset of $A^\circ$. Nets (here with parameter $\e$) are written in the form $(u_\e)_\e$. The class with respect to any equivalence relation is denoted by square brackets $[\dotsc]$. A family of objects $x_i$ indexed by $i \in I$ is written as $\{x_i\}_{i \in I}$ or simply $\{ x_i \}_i$ when the index set is clear from the context. We use Landau notation: for expressions $f(\e)$ and $g(\e)$ depending on and defined for small $\e$ we write $f(\e) = O(g(\e))$ (always for $\e \rightarrow 0)$ if and only if $\exists C>0$ $\exists \e_0>0$ $\forall \e<\e_0$: $\abso{f(\e)} < C \abso{g(\e)}$. $B_\eta(x)$ resp.~ $B_\eta(K)$  denotes the metric ball of radius $\eta$ around $x \in \bR^n$ resp.~ a set $K$, $\dist$ denotes the Euclidean distance function on $\bR^n$. For a function $f(\varphi,x)$ of a variable $\varphi$ and an $n$-dimensional real variable $x=(x_1,\dotsc,x_n)$, $\ud_2 f$ denotes the total differential of $f$ with respect to $x$ and $\pd_i f$ its partial differential with respect to $x_i$. For the derivative of a function $\gamma$ depending on $t \in \bR$ we will write $\gamma'$. An $n$-tuple $\alpha=(\alpha_1,\dotsc, \alpha_n) \in \bN_0^n$ is called a multi-index; we use the notation $\abso{\alpha} = \alpha_1+\dotsc+\alpha_n$, $x^\alpha = x_1^{\alpha_1}\dotsm x_n^{\alpha_n}$, and $\pd^\alpha f = \pd_1^{\alpha_1} \dotsc \pd_n^{\alpha_n}$. A strictly decreasing sequence $(x_n)_{n \in \bN}$ converging to $x_0$ is denoted by $x_n \searrow x_0$. A function between finite dimensional real vector spaces is said to be smooth if it is infinitely differentiable. The action of a distribution $u \in \ccD'(\Omega)$ on a test function $\varphi \in \ccD(\Omega)$ is written as $\langle u, \varphi \rangle$.

\subsection{Calculus on convenient vector spaces}

The construction of the diffeomorphism invariant full algebra $\cG^d(\Omega)$ as defined below requires calculus on infinite-dimensional locally convex spaces as an indispensable prerequisite. The theoretical framework chosen for this by Grosser et al.~\cite{GKOS} is smooth calculus on convenient vector spaces, which is presented by Kriegl and Michor \cite{KM} using functional analysis and by Fr\"olicher and Kriegl \cite{Froe} using category theory. For a detailed exposition of what is needed for the diffeomorphism invariant full algebra we refer to \cite{GKOS}, Section 2.2. Whenever we encounter smoothness on a subset of a locally convex space (or an affine subspace thereof) we endow it it with the initial smooth structure.

A sesquilinear form on a complex locally convex space is smooth if and only if it is bounded; this easily results from an adaptation of \cite{KM} Section 5 to antilinear maps.

Although the differential is at first only defined for mappings having as domain open subsets of locally convex spaces with respect to a certain topology (\cite{KM} Theorem 3.18) this definition can be easily extended to maps defined on affine subspaces, as is remarked in the proof of Proposition \ref{srxprop}. Properties like the chain rule and the symmetry of higher derivatives remain intact.

\subsection{Colombeau Algebras}

We will now give the definitions of the special algebra $\cG^s(\Omega)$ and the full algebras $\cG^e(\Omega)$ and $\cG^d(\Omega)$ on an arbitrary open subset $\Omega \subseteq \bR^n$.

The special Colombeau algebra $\cG^s(\Omega)$ (\cite{GKOS} Section 1.2) consists of nets of smooth functions on $\Omega$ indexed by $I \coleq (0,1]$. Such a net $(u_\e)_\e \in C^\infty(\Omega)^I$ is said to be \textit{moderate} if $\forall K \csub \Omega$ $\forall \alpha \in \bN_0^n$ $\exists N \in \bN$ such that $\sup_{x \in K}\abso{\pd^\alpha u_\e(x)} = O(\e^{-N})$, or \textit{negligible} if $\forall K \csub \Omega$ $\forall \alpha \in \bN_0^n$ $\forall m \in \bN:$ $\sup_{x \in K}\abso{\pd^\alpha u_\e(x)} = O(\e^m)$. $\cG^s(\Omega)$ then is the quotient of $\cE^s_M(\Omega)$ (the set of moderate nets) modulo $\cN^s(\Omega)$ (the set of negligible nets).

The full algebras $\cG^e(\Omega)$ and $\cG^d(\Omega)$ require some auxiliary definitions. For $q \in \bN_0$ let $\cA_q(\Omega)$ be the set of all test functions $\varphi \in \ccD(\Omega)$ having integral 1, if $q \ge 1$ additionally satisfying $\int x^\alpha \varphi(x)\,\ud x = 0$ for all multi-indices $\alpha$ with $1 \le \abso{\alpha} \le q$. Let $\cA_{0q}(\Omega)$ be defined in the same way but with integral $0$. For any subset $M \subseteq \Omega$ define $\cA_{0,M}(\Omega)$ as the set of those elements of $\cA_0(\Omega)$ with support in $M$. $\cA_q(\Omega)$ and $\cA_{0q}(\Omega)$ are endowed with the initial topology and the initial smooth structure with respect to the embedding in $\ccD(\Omega)$ or $\ccD(\bR^n)$. Let $U(\Omega)$ be the set of all pairs $(\varphi, x) \in \cA_0(\bR^n) \times \Omega$ satisfying $\supp \varphi +x \subseteq \Omega$. Furthermore, let $C^\infty_b(I \times \Omega, \cA_0(\bR^n))$ be the space of those mappings which are smooth from $I \times \Omega$ into $\cA_0(\bR^n)$ such that for any compact set $K \csub \Omega$ and any $\alpha \in \bN_0^n$ the set $\{ \pd^\alpha \phi(\e,x)\ |\ \e \in I,\ x \in K\}$ is bounded in $\ccD(\bR^n)$. For $\e \in \bR_+$ let $S_\e: \ccD(\bR^n) \rightarrow \ccD(\bR^n)$ be the mapping given by $(S_\e\varphi)(y) \coleq \e^{-n}\varphi(y/\e)$ and set $S^{(\e)}(\varphi, x) \coleq (S_\e\varphi, x)$ for $(\varphi, x) \in \ccD(\bR^n) \times \bR^n$. For $x \in \bR^n$ denote by $T_x: \ccD(\bR^n) \rightarrow \ccD(\bR^n)$ the mapping given by $(T_x\varphi)(y) \coleq \varphi(y-x)$ and define $T: \ccD(\Omega) \times \bR^n \rightarrow \ccD(\Omega) \times \bR^n$ by $T(\varphi, x) \coleq (T_x\varphi, x)$. For a map $R$ we will frequently write $R_\e$ instead of $R \circ S^{(\e)}$.

For $\cG^e(\Omega)$ (\cite{GKOS} Section 1.4), the base space $\cE^e(\Omega)$ is the set of all functions $R: U(\Omega) \rightarrow \bC$ which are smooth in the second variable. $R$ is called \textit{moderate} if $\forall K \csub \Omega$ $\forall \alpha \in \bN_0^n$ $\exists N \in \bN$ $\forall \varphi \in \cA_N(\bR^n)$: $\sup_{x \in K}\abso{\pd^\alpha R(S_\e\varphi, x)} = O(\e^{-N})$ and \textit{negligible} if $\forall K \csub \Omega$ $\forall \alpha \in \bN_0^n$ $\forall m \in \bN$ $\exists q \in \bN$ $\forall \varphi \in \cA_q(\bR^n)$: $\sup_{x \in K}\abso{\pd^\alpha R(S_\e\varphi, x)} = O(\e^m)$. The corresponding sets $\cE^e_M(\Omega)$ of moderate and $\cN^e(\Omega)$ of negligible functions give rise to the differential algebra $\cG^e(\Omega) \coleq \cE^e_M(\Omega) / \cN^e(\Omega)$. Distributions $u \in \ccD'(\Omega)$ are embedded via the linear injective mapping $\iota:\ccD'(\Omega) \rightarrow \cE^e_M(\Omega)$ given by $\iota(u)(\varphi,x) \coleq \langle u, T_x\varphi \rangle$ for $(\varphi,x) \in U(\Omega)$. The derivations of $\cG^e(\Omega)$ which extend the distributional ones are given by $(D_iR)(\varphi, x) \coleq (\pd_i R)(\varphi, x)$ for $R \in \cE^e_M(\Omega)$ and $i=1,\dotsc,n$.

For $\cG^d(\Omega)$ (\cite{GKOS} Chapter 2), we take the base space $\cE^d(\Omega) \coleq C^\infty(U(\Omega))$. A map $R \in \cE^d(\Omega)$ is called \textit{moderate} if $\forall K \csub \Omega$ $\forall \alpha \in \bN_0^n$ $\exists N \in \bN$ $\forall \phi \in C^\infty_b(I \times \Omega, \cA_0(\bR^n))$: $\sup_{x \in K}\abso{\pd^\alpha R(S_\e\phi(\e,x), x)} = O(\e^{-N})$ and \textit{negligible} if it is moderate and $\forall K \csub \Omega$ $\forall \alpha \in \bN_0^n$ $\forall m \in \bN$ $\exists q \in \bN$ $\forall \phi \in C^\infty_b(I \times \Omega, \cA_q(\bR^n))$: $\sup_{x \in K}\abso{\pd^\alpha R(S_\e\phi(\e,x), x)} = O(\e^m)$. The corresponding sets $\cE^d_M(\Omega)$ of moderate and $\cN^d(\Omega)$ of negligible functions give rise to the differential algebra $\cG^d(\Omega) \coleq \cE^d_M(\Omega) / \cN^d(\Omega)$. The embedding (denoted by $\iota$ as well) of distributions $u \in \ccD'(\Omega)$ is given by $\iota(u)(\varphi,x) \coleq \langle u, T_x\varphi \rangle$ for $(\varphi,x) \in U(\Omega)$. The derivations which extend the distributional ones are given by $(D_iR)(\varphi, x) \coleq (\pd_iR)(\varphi,x)$.

A \textit{constant} in one of the preceding differential algebras (as in any differential ring) is an element whose derivations are all zero (\cite{Kolchin} Chapter I Section 1).

\begin{remark}\label{remarkuno}
\begin{enumerate}
\item[(i)]In all definitions of moderateness and negligibility above and below, when expanding the Landau Symbol in expressions of the form $\abso{f(\e)} = O(\e^{-N})$ into $\exists C>0$ $\exists \eta>0$ $\forall \e<\eta$: $\abso{f(\e)} < C\e^{-N}$ (resp.~ $\e^m$ for negligibility) one can without loss of generality fix $C=1$ without changing the definitions.
\item[(ii)] In the definitions of negligibility one can disregard the derivatives and only consider $\alpha=0$ if one presupposes the tested element to be moderate (\cite{GKOS} Theorems 1.2.3, 1.4.8, and 2.5.4).
\end{enumerate}
\end{remark}

\section{Previous results in the special algebra $\cG^s(\Omega)$}

We first recall the definition of generalized points, numbers, and point values for $\cG^s(\Omega)$. Two results justify these definitions: first, the ring of constants in $\cG^s(\Omega)$ equals the space of generalized numbers. Second, two generalized functions are equal if and only if they have the same point values.

\begin{definition}[\cite{GKOS} Definition 1.2.31] Generalized numbers in the $\cG^s$-setting are defined by
\begin{align*}
\bC_M &\coleq \{ (r_\e)_\e \in \bC^I\ |\ \exists N \in \bN:\ \abso{r_\e} = O(\e^{-N}) \},\\
\bC_N &\coleq \{ (r_\e)_\e \in \bC^I\ |\ \forall m \in \bN:\ \abso{r_\e} = O(\e^m) \},\\
\widetilde\bC &\coleq \bC_M / \bC_N.
\end{align*}
\end{definition}

\begin{definition}[\cite{GKOS} Definition 1.2.44] Generalized points in the $\cG^s$-setting are defined by
\begin{align*}
\Omega_M & \coleq \{ (x_\e)_\e \in \Omega^I\ |\ \exists N \in \bN:\ \abso{x_\e} = O(\e^{-N}) \},\\
(x_\e)_\e \sim (y_\e)_\e &\coliff \forall m \in \bN:\ \abso{x_\e - y_\e} = O(\e^m),\\
\widetilde \Omega & \coleq \Omega_M / \sim,\\
\widetilde\Omega_c & \coleq \{ \tilde x = [(x_\e)_\e] \in \widetilde\Omega\ |\ \exists K \csub \Omega\ \exists \eta>0\ \forall \e<\eta:\ x_\e \in K \}.
\end{align*}
\end{definition}

Clearly $\widetilde \bC$ can be seen as a subset of $\cG^s(\Omega)$.

\begin{proposition}[\cite{GKOS} Proposition 1.2.35]Let $\Omega \subseteq \bR^n$ be connected and $\tilde u \in \cG^s(\Omega)$. Then $D\tilde u = 0$ if and only if  $\tilde u \in \widetilde\bC$. 
\end{proposition}

\begin{definition}Let $\tilde u=[(u_\e)_\e] \in \cG^s(\Omega)$ and $\tilde x = [(x_\e)_\e] \in \widetilde \Omega_c$. Then the generalized point value of $\tilde u$ at $\tilde x$ defined by $\tilde u(\tilde x) \coleq [(u_\e(x_\e))_\e]$ is a well-defined element of $\widetilde\bC$.
\end{definition}

\begin{theorem}[\cite{GKOS} Theorem 1.2.64]Let $\tilde u \in \cG^s(\Omega)$. Then $\tilde u = 0$  in $\cG^s(\Omega)$ if and only if $\tilde u(\tilde x) = 0 \textrm{ in }\widetilde\bC$ for all $\tilde x \in \widetilde \Omega_c$.
\end{theorem}

\section{Point values in $\cG^e(\Omega)$}

It was asserted by Grosser et al.~(\cite{GKOS} Section 1.4.2) that results concerning point values obtained in the special algebra can be recovered in the full algebra $\cG^e(\Omega)$. This section explicitly states these results and their proofs for $\cG^e(\Omega)$, which should not be regarded as a mere technical exercise but as an essential building step if one aims to obtain the corresponding results in $\cG^d(\Omega)$, where in addition one needs to incorporate smoothness into the proofs presented here.

After recalling the definition of generalized numbers in the $\cG^e$-setting we will define a suitable space of generalized points.

\begin{definition}[\cite{GKOS} Definition 1.4.19]\label{genzge}Generalized numbers in the $\cG^e$-setting are defined by
\begin{align*}
\bC_M(n) &\coleq \{ r: \cA_0(\bR^n) \rightarrow \bC\ |\ \exists N \in \bN\ \forall \varphi \in \cA_N(\bR^n):\\
&\qquad\qquad \abso{r(S_\e\varphi)} = O(\e^{-N}) \},\\
\bC_N(n) &\coleq \{ r: \cA_0(\bR^n) \rightarrow \bC\ |\ \forall m \in \bN\ \exists q \in \bN\ \forall \varphi \in \cA_q(\bR^n):\\
&\qquad\qquad \abso{r(S_\e\varphi)} = O(\e^m) \},\\
\widetilde\bC(n) &\coleq \bC_M(n) / \bC_N(n).
\end{align*}
\end{definition}

\begin{definition}Generalized points in the $\cG^e$-setting are defined by
\begin{align*}
\begin{split}
\Omega_M(n) &\coleq \{ X: \cA_0(\bR^n) \rightarrow \Omega\ |\ \forall \varphi \in \cA_0(\bR^n)\ \exists \e_0>0\ \forall \e<\e_0:\\
&\qquad (S_\e\varphi, X(S_\e\varphi)) \in U(\Omega)\text{ and}\\
&\qquad \exists N\in \bN\ \forall \varphi \in \cA_N(\bR^n):\ \abso{X(S_\e\varphi)} = O(\e^{-N})\},
\end{split}\\
\Omega_N(n) &\coleq \{ X \in \Omega_M(n)\ |\ \forall m\in \bN\ \exists q\in \bN\ \forall \varphi\in\cA_q(\bR^n): \\
&\qquad \abso{X(S_\e\varphi)} = O(\e^m),\\
\widetilde \Omega(n) &\coleq \Omega_M(n) / \Omega_N(n),\\
\begin{split}
  \widetilde \Omega_c(n) &\coleq \{ \widetilde X \in \widetilde\Omega(n)\ |\ \text{ for one (thus any) representative }X\text{ of }\widetilde X\\
&\qquad \exists K \csub \Omega\ \exists N \in \bN\ \forall \varphi \in \cA_N(\bR^n)\\
&\qquad \exists \eta>0\ \forall \e<\eta:\ X(S_\e\varphi) \in K \}.
\end{split}
\end{align*}
We write $X \sim Y$ if $X-Y \in \Omega_N(n)$. Any $X \in \Omega_M(n)$ satisfying the condition in the definition of $\widetilde \Omega_c(n)$ is called \textit{compactly supported} (in $K$).
If one replaces $\bC$ by $\bR$ in Definition \ref{genzge} the resulting space is denoted by $\widetilde\bR(n)$.

\end{definition}

\begin{proposition}\label{gepvwd}Let $X \in \Omega_M(n)$ be compactly supported and let $R \in \cE_M^e(\Omega)$. Define $R(X): \cA_0(\bR^n) \rightarrow \bC$ by
\begin{equation*}
R(X)(\varphi) \coleq \left\{ \begin{aligned}
R(\varphi,X(\varphi)) & \qquad (\varphi, X(\varphi)) \in U(\Omega)\\
0 & \qquad \textrm{otherwise}.
\end{aligned}\right.
\end{equation*}
Then $R(X)$ is in $\bC_M(n)$, $R \in \cN^e(\Omega)$ implies $R(X) \in \bC_N(n)$, and $X \sim Y$ implies $R(X) - R(Y) \in \bC_N(n)$.
\end{proposition}
\begin{proof}
Let $X$ be compactly supported in $K \csub \Omega$, which means that $\exists N \in \bN$ $\forall \varphi \in \cA_{N}(\bR^n)$: $X(S_\e\varphi) \in K$ for small $\e$. Given any $\varphi \in \cA_N(\bR^n)$, for small $\e$ we have $X(S_\e\varphi) \in K$, $R(X)(S_\e\varphi) = R(S_\e\varphi, X(S_\e\varphi))$, and thus $\abso{R(X)(S_\e\varphi)} \le \sup_{x \in K}\abso{R(S_\e\varphi, x)}$ whence $R(X)$ inherits moderateness respectively negligibility from $R$. For the last claim, choose some $m \in \bN$ for the test for membership in $\bC_N(n)$. Then we use the following ingredients.
\begin{enumerate}
\item[(i)] As $X \sim Y$, $\exists q_0 \in \bN$ $\forall \varphi \in \cA_{q_0}(\bR^n)$: $\abso{X(S_\e\varphi)-Y(S_\e\varphi)} < \e^m$ for small $\e$.
\item[(ii)] $\exists \eta>0$: $\overline{B_\eta(K)}\subseteq \Omega$. Set $V \coleq B_\eta(K)$.
\item[(iii)] As derivatives of $R$ are moderate, $\exists N' \in \bN$ $\forall \varphi \in \cA_{N'}(\bR^n)$ such that $\sup_{x \in \overline{V}}\abso{\ud_2 R(S_\e\varphi, x)} \le \e^{-N'}$ for small $\e$.
\item[(iv)] From (i) we know in particular that given $\varphi \in \cA_{\max(q_0, N)}(\bR^n)$, $g(t) \coleq (X + t(Y - X))(S_\e\varphi)$ lies in $V$ for small $\e$ and all $t \in [0,1]$.
\item[(v)] $\forall \varphi \in \cA_0(\bR^n)$: $\supp S_\e\varphi + V \subseteq \Omega$ for small $\e$.
\end{enumerate}

Next let $\varphi \in \cA_{\max(q_0,N,N')}(\bR^n)$ and $\e$ small enough. Then by (iv), $X(S_\e\varphi)$ and $Y(S_\e\varphi)$ are in $V$, $(R(X) - R(Y))(S_\e\varphi) = R(S_\e\varphi,X(S_\e\varphi)) - R(S_\e\varphi, Y(S_\e\varphi))$ and the domain of $R(S_\e\varphi, \cdot)$ contains $V$. Set $F(t) \coleq R(S_\e\varphi, g(t))$ for $t \in [0,1]$. Then $F$ is smooth on $[0,1]$ and
\begin{multline*}
\abso{R(X)(S_\e\varphi) - R(Y)(S_\e\varphi)} = \abso{F(1)-F(0)} = \abso{\int_0^1 F'(t)\ud t} \\
= \abso{\int_0^1 \ud_2 R(S_\e\varphi, g(t)) \cdot (X(S_\e\varphi)-Y(S_\e\varphi))\ud t} \\
\le \abso{(X-Y)(S_\e\varphi)} \cdot \sup_{x \in \overline{V}}\abso{(\ud_2 R)(S_\e\varphi, x)} \leq \e^m \e^{-N'}.
\end{multline*}
As $m$ was arbitrary this concludes the proof.
\end{proof}

The following lemma will be used to construct generalized points and numbers taking prescribed values.

\begin{lemma}\label{constructpoint}Given $\varphi_q \in \cA_q(\bR^n)$, $\e_{q,k} \in (0,\infty)$ and $x_0$, $x_{q,k}$ in any set $A$ for all $q,k \in \bN$, there exists a mapping $X: \cA_0(\bR^n) \rightarrow A$ and strictly increasing sequences $(q_l)_{l \in \bN}$ and $(a_l)_{l \in \bN}$ of natural numbers such that $X(S_{\e_{q_l, k}} \varphi_{q_l}) = x_{q_l,k}$ $\forall k,l \in \bN$, $X(\varphi) = x_0$ for all $\varphi$ not equal to some $S_{\e_{q_l, k}}\varphi_{q_l}$, and $\varphi_{q_l} \in \cA_{a_l}(\bR^n) \setminus \cA_{a_l+1}(\bR^n)$.
\end{lemma}
\begin{proof}
  Set $q_1 \coleq 1$, $a_1$ such that $\varphi_{q_1} \in \cA_{a_1}(\bR^n) \setminus \cA_{a_1+1}(\bR^n)$ and inductively choose $q_{l+1} \coleq a_l+1$ and $a_{l+1}$ appropriately. This is possible because for $q$ increasing more and more moments of $\varphi_q$ have to vanish. Then define $X: \cA_0(\bR^n) \rightarrow A$ as follows: given $\psi \in \cA_0(\bR^n)$, if $\psi = S_{\e_{q_l,k}}\varphi_{q_l}$ for some $k,l$ then set $X(\psi) \coleq x_{q_l, k}$, otherwise set $X(\psi) \coleq x_0$.
\end{proof}

\begin{definition}For $\widetilde R = [R] \in \cG^e(\Omega)$ and $\widetilde X = [X] \in \widetilde\Omega_c(n)$ we define the point value $\widetilde R(\widetilde X)$ of $\widetilde R$ at $\widetilde X$ as the class in $\widetilde\bC(n)$ of $R(X)$ as defined in Proposition \ref{gepvwd}.
\end{definition}

Having defined suitable spaces of generalized points and numbers as well as a corresponding notion of point evaluation we can now state the point value characterization theorem for $\cG^e$.

\begin{theorem}\label{gepvchar}Let $\widetilde R = [R] \in \cG^e(\Omega)$. Then $\widetilde R=0$ if and only if $\widetilde R (\widetilde X)=0$ in $\widetilde\bC(n)$ for all $\widetilde X \in \widetilde\Omega_c(n)$.
\end{theorem}
\begin{proof}Necessity was already shown in Proposition \ref{gepvwd}. For sufficiency assume that $R \not\in \cN^e(\Omega)$; then by Remark \ref{remarkuno} (ii) there exist $K \csub \Omega$ and $m_0 \in \bN$ such that for all $q \in \bN$ there is some $\varphi_q \in \cA_q(\bR^n)$, a sequence $(\e_{q,k})_{k \in \bN} \searrow 0$ and a sequence $(x_{q,k})_{k \in \bN}$ in $K$ such that $\abso{R(S_{\e_{q,k}}\varphi_q, x_{q,k})} \ge \e_{q,k}^{m_0}$.

  Let $X:\cA_0(\bR^n) \rightarrow K$, $(q_l)_{l \in \bN}$ and $(a_l)_{l \in \bN}$ be as obtained from Lemma \ref{constructpoint} with arbitrary $x_0 \in K$.
 Then clearly $X$ is compactly supported, $[X] \in \widetilde\Omega_c$ and $R(X) \not\in \bC_N(n)$: for any $q \in \bN$ there is some $l \in \bN$ such that $a_l\ge q$, so $\varphi_{q_l} \in \cA_q(\bR^n)$. By construction,
\begin{align*}
\abso{R(X)(S_{\e_{q_l, k}}\varphi_{q_l})} &= \abso{R(S_{\e_{q_l, k}} \varphi_{q_l}, X(S_{\e_{q_l, k}}\varphi_{q_l}))} \\
& = \abso{R(S_{\e_{q_l, k}} \varphi_{q_l}, x_{q_l, k})} > \e^{m_0}_{q_l, k}
\end{align*}
for all large enough $k \in \bN$, which ensures that the negligibility test for $R(X)$ fails.
\end{proof}

The proof of the following Proposition is evident.

\begin{proposition}The map $\rho: \bC_M(n) \rightarrow \cE^e(\Omega)$ given by $\rho(r)(\varphi, x) \coleq r(\varphi)$ $\forall (\varphi, x) \in U(\Omega)$ is a ring homomorphism preserving moderateness and negligibility and thus induces an embedding $\tilde \rho: \widetilde \bC(n) \rightarrow \cG^e(\Omega)$.
\end{proposition}

\begin{lemma}\label{compblah}
Let $\Omega \subseteq \bR^n$ be connected and $K \csub \Omega$. Then there exist a set $M \csub \Omega$ containing $K$ and a real number $L>0$ such that any two points in $K$ can be connected by a continuous curve $\gamma: [0,1] \rightarrow \Omega$
with image in $M$ having length $\int_0^1 \abso{ \gamma'(t)} \ud t \leq L$.
\end{lemma}
\begin{proof}
Cover $K$ by finitely many closed balls of some radius $\e>0$ which are contained in $\Omega$. As $\Omega$ is (pathwise) connected they can be joined by finitely many continuous curves in $\Omega$. Taking as $M$ the union of these $\e$-balls and the images of these curves, the existence of $L$ as desired is obvious.
\end{proof}

In the differential algebra $\cG^e(\Omega)$ the constant elements are exactly those whose derivatives are zero. With the availability of point values one can also call a generalized function constant if it has the same generalized value at every generalized point. The following proposition shows that these properties in fact are equivalent.

\begin{proposition}If $\widetilde R \in \cG^e(\Omega)$ has the property $\widetilde R (\widetilde X) = \widetilde R(\widetilde Y)$ $\forall \widetilde X,\widetilde Y \in \widetilde\Omega_c(n)$ then $D_i\widetilde R=0$ for $i=1,\dotsc,n$; if $\Omega$ is connected the converse also holds.
\end{proposition}
\begin{proof}
Given any $\widetilde X \in \widetilde\Omega_c(n)$, one easily sees that for all $\widetilde Y \in \widetilde\Omega_c(n)$ we have $\tilde \rho(\widetilde R(\widetilde X))(\widetilde Y) = \widetilde R(\widetilde X)$ on the one hand, and $\widetilde R(\widetilde Y)=\widetilde R(\widetilde X)$ on the other hand by assumption. By Theorem \ref{gepvchar} then $\tilde\rho(\widetilde R(\widetilde X))=\widetilde R$, whence $D_i\widetilde R = D_i\tilde\rho(\widetilde R(\widetilde X)) = 0$ follows at once from the definitions.

For the converse we show that in case $\Omega$ is connected $D_i\widetilde R = 0$ (for $i=1,\dotsc,n)$ in $\cG^e(\Omega)$ implies $\widetilde R = \tilde\rho(\widetilde R(\widetilde X))$ for arbitrary $\widetilde X = [X] \in \widetilde\Omega_c(n)$. Fix $K_1 \csub \Omega$ and $m \in \bN$ for testing and let $X$ be compactly supported in $K_2 \csub \Omega$. Let $M$ and $L$ be as obtained from Lemma \ref{compblah} applied to $K = K_1 \cup K_2$. By assumption,
\begin{enumerate}
  \item[(i)] $\exists q \in \bN$ $\forall \varphi \in \cA_q(\bR^n)$ $\exists \e_0>0$ $\forall \e<\e_0$: $\sup_{x \in M} \abso{ \ud_2 R(S_\e\varphi, x)} \leq \e^m$.
\item[(ii)] $\exists N \in \bN$ $\forall \varphi \in \cA_N(\bR^n)$ $\exists \eta>0$ $\forall \e < \eta$: $X(S_\e\varphi) \in K_2$.
\end{enumerate}
Now let $\varphi \in \cA_{\max(q,N)}(\bR^n)$ and $\e < \min(\e_0,\eta)$. Then for every $y \in K_1$ there exists a continuous curve $\gamma: [0,1] \rightarrow \Omega$ with image in $M$ connecting $y$ and $X(S_\e\varphi)$ and having length $\le L$. Thus we can estimate 
\begin{align*}
  \abso{R(S_\e\varphi, y)-R(S_\e\varphi, X(S_\e\varphi))} &= \abso{\int_0^1 \ud_2 R(S_\e\varphi, \gamma(t))\gamma'(t)\ud t}\\
  & \leq \sup_{x \in M} \abso{ \ud_2 R(S_\e\varphi, x) } \cdot \int_0^1 \abso{\gamma'(t)} \ud t \leq L\e^m
\end{align*}
which gives the claimed result.
\end{proof}

\begin{definition}For $\tilde r,\tilde s \in \widetilde\bR(n)$ we write $\tilde r \le \tilde s$ if there are representatives $r$, $s$ such that $r(\varphi) \le s(\varphi)$ for all $\varphi \in \cA_0(\Omega)$.
\end{definition}

\begin{proposition} $(\widetilde\bR(n), \leq)$ is a partially ordered ring.
\end{proposition}

\begin{proof}Reflexivity is clear. For antisymmetry, $\tilde r \leq \tilde s$ and $\tilde s\leq \tilde r$ imply $r_1 \leq s_1$ and $s_2 \leq r_2$ for some representatives $r_1,r_2$ of $\tilde r$ and $s_1,s_2$ of $\tilde s$. Writing $s_1=s_2 + n$ and $r_2=r_1+m$ with $n,m \in \cN^e(\Omega)$ gives $r_1 \leq s_2+n$ and $s_2 \leq r_1+m$, thus $r_1-s_2 \leq n$ and $s_2 - r_1 \leq m$, implying $\abso{r_1-s_2} \leq \max(n,m)$ and finally $r_1-s_2 \in \cN^e(\Omega)$. For transitivity assume $\tilde r \leq \tilde s \leq \tilde t$. Then with representatives $r, s_1, s_2$ and $t$ we get $s_1 = s_2+n$ with $n \in \cN^e(\Omega)$ and thus $r \leq s_1 = s_2+n \leq t+n$, which is $\tilde r \le \tilde t$. Finally, $\tilde r \leq \tilde s$ clearly implies $\tilde r + \tilde t \leq \tilde s + \tilde t$ and $0 \leq \tilde r$, $0 \leq \tilde s$ reads $n \leq r$, $m \leq s$ in representatives which implies $nm \leq rs$ or $0 \le \tilde r \tilde s$.
\end{proof}

We call a generalized number $\tilde r \in \widetilde\bC(n)$ strictly nonzero if it has a representative $r \in \bC_M(n)$ such that
\begin{equation}\label{strictlynonzero}
\exists q\in \bN\ \forall \varphi \in \cA_q(\bR^n)\ \exists C>0\ \exists \eta>0\ \forall \e<\eta:\ \abso{r(S_\e\varphi)} > C\e^q.
\end{equation}

Note that Remark \ref{remarkuno} (i) applies here and we can always have $C=1$. We get the following characterization of invertibility in $\widetilde \bC(n)$.

\begin{proposition}\label{inverti}An element of $\widetilde\bC(n)$ is invertible if and only if it is strictly nonzero.
\end{proposition}

\begin{proof}
Given $\tilde r = [r],\tilde s=[s] \in \widetilde\bC(n)$ with $\tilde r\tilde s=1$, there exists $t \in \bC_N(n)$ such that $rs = 1 + t$. By the definition of negligibility $\exists q\in\bN$ $\forall \varphi \in \cA_q(\bR^n)$ $\exists \eta>0$ $\forall \e<\eta$: $\abso{t(S_\e\varphi)} < 1/2$, and thus also $s(S_\e\varphi) \ne 0$.
By moderateness of $s$ $\exists N \in \bN$ $\forall \varphi \in \cA_N(\bR^n)$ $\exists \eta'>0$ such that for all $\e < \eta'$ we have $\abso{s(S_\e\varphi)} < \e^{-N}$.
Thus for $q' \coleq \max(q,N)$, $\varphi \in \cA_{q'}(\bR^n)$, and $\e < \min(\eta, \eta')$ we obtain
\[ \abso{r(S_\e\varphi)} = \abso{\frac{1 + t(S_\e\varphi)}{s(S_\e\varphi)}} > \frac{\e^N}{2}\ge \frac{\e^{q'}}{2}. \]

Conversely, given $r \in \bC_M(n)$ satisfying \eqref{strictlynonzero} set $s(\varphi) \coleq 1/r(\varphi)$ where defined and $0$ elsewhere. Then $s \in \bC_M(n)$ by definition and obviously $rs-1 \in \cN^e(n)$ because for $\varphi \in \cA_q(\bR^n)$ with $q$ of \eqref{strictlynonzero} and small $\e$, $s(S_\e\varphi)=1/r(S_\e\varphi)$, thus $rs-1=0$ and the negligibility test succeeds trivially.
\end{proof}

\begin{proposition}For $\tilde r \in \widetilde\bC(n)$ the following assertions are equivalent.
\begin{enumerate}
\item[(i)] $\tilde r$ is not invertible.
\item[(ii)] $\tilde r$ has a representative $r$ such that for all $q \in \bN$ there is some $\varphi_q \in \cA_q(\bR^n)$ and a sequence $(\e_{q,k})_{k \in \bN} \searrow 0$ such that $r(S_{\e_{q,k}}\varphi_q)=0$ for all $k$.
\item[(iii)] $\tilde r$ is a zero divisor.
\end{enumerate}
\end{proposition}

\begin{proof}
(i) $\Rightarrow$ (ii): $\tilde r$ fails to be strictly nonzero, thus any representative $r$ satisfies
\[ \forall q \in \bN\ \exists \varphi_q \in \cA_q(\bR^n)\ \exists (\e_{q,k})_{k \in \bN} \searrow 0:\ \abso{r(S_{\e_{q,k}}\varphi_q)} \le \e_{q,k}^q. \]
With $x_{q,k} \coleq r(S_{\e_{q,k}}\varphi_q)$ ($q,k \in \bN$) and $x_0 \coleq 0$ let $s: \cA_0(\bR^n) \rightarrow \bC$, $(q_l)_l$, and $(a_l)_l$ be as obtained from Lemma \ref{constructpoint}. This map $s$ satisfies $s(S_{\e_{q_l, k}}\varphi_{q_l}) = x_{q_l, k}$ $\forall k,l \in \bN$. Then $s$ is negligible: let $m \in \bN$ be given and choose $l_0 \in \bN$ such that $q_{l_0} > m$. Let $\varphi \in \cA_{a_l}(\bR^n)$. Then $s(S_\e\varphi)$ can only be nonzero if $\varphi = S_\eta \varphi_{q_l}$ for some $\eta>0$ and $l \ge l_0$ and this requires that $S_\e\varphi = S_\e S_\eta \varphi_{q_l} = S_{\e_{q_l,k}}\varphi_{q_l}$ for some $k \in \bN$, that is $\e\eta = \e_{q_l, k}$. In this case
\[ \abso{s(S_\e\varphi)} = \abso{r(S_{\e_{q_l, k}}\varphi_{q_l})} \le \e_{q_l,k}^{q_l} = \eta^{q_l}\e^{q_l} < \eta^{q_l}\e^m \]
for all $\e = \e_{q_l, k}/\eta$ which are $<1$. Finally $r-s$ has the desired property: given $q \in \bN$, there is some $l$ such that $q_l \ge q$ and for $\varphi_{q_l} \in \cA_{q_l}(\bR^n) \subseteq \cA_q(\bR^n)$ we have $(r-s)(S_{\e_{q_l,k}}\varphi_{q_l}) = 0$.

(ii) $\Rightarrow$ (iii): Define $s: \cA_0(\bR^n) \rightarrow \bC$ by $s(\varphi) \coleq 1$ if $r(\varphi)=0$ and $s(\varphi)\coleq 0$ otherwise. Then $s \in \bC_M(n)$ and $rs=0$ but it is easily verified that $s \not\in \bC_N(n)$.

(iii) $\Rightarrow$ (i) is trivial.
\end{proof}

The following is a characterization of non-degeneracy of matrices over $\widetilde\bC(n)$.

\begin{proposition}Let $A \in \widetilde\bC(n)^{m^2}$ be an $m \times m$ square matrix with entries from $\widetilde\bC(n)$. The following are equivalent:
\begin{enumerate}
\item[(i)] $A$ is non-degenerate, i.e., if $\xi,\eta \in \widetilde\bC(n)^m$ then $\xi^tA\eta = 0$ $\forall \eta$ implies $\xi=0$.
\item[(ii)] $A: \widetilde\bC(n)^m \rightarrow \widetilde\bC(n)^m$ is injective.
\item[(iii)] $A: \widetilde\bC(n)^m \rightarrow \widetilde\bC(n)^m$ is bijective.
\item[(iv)] $\det(A)$ is invertible.
\end{enumerate}
\end{proposition}
\begin{proof}The proof is purely algebraical and hence is entirely equivalent to the version for $\cG^s(\Omega)$ (\cite{GKOS} Lemma 1.4.41). More explicitly, (ii) $\Leftrightarrow$ (iii) $\Leftrightarrow$ (iv) is dealt with by \cite{Bourbaki} Chapter III \S 8 Proposition 3 and Theorem 1. (i) $\Rightarrow$ (ii) follows by showing that (i) is equivalent to $A^t$ being injective, after which (ii) $\Rightarrow$ (iv) can be applied to $\det(A) = \det(A^t)$.
\end{proof}

The next theorem is a characterization of invertibility of generalized functions in $\cG^e(\Omega)$.

\begin{theorem}\label{geinv}For $\widetilde R \in \cG^e(\Omega)$ the following are equivalent:
\begin{itemize}
\item[(i)] $\widetilde R$ is invertible.
\item[(ii)] For each representative $R$ of $\widetilde R$ the following holds:
\begin{multline*}
\forall K \csub \Omega\ \exists m \in \bN\ \exists q \in \bN\ \forall \varphi \in \cA_q(\bR^n)\\
 \exists C>0\ \exists \e_0>0\ \forall \e<\e_0: \sup_{x \in K} \abso{R(S_\e\varphi, x)} > C \e^m. 
\end{multline*}
Remark \ref{remarkuno} (i) applies here; furthermore, we can always have $m=q$.
\end{itemize}
\end{theorem}
\begin{proof}
Assuming (i), there exist $S \in \cE^e_M(\Omega)$ and $Q \in \cE^e_N(\Omega)$ such that $RS = 1+Q$. Fix $K \csub \Omega$. Then $\exists p \in \bN$ $\forall \varphi \in \cA_p(\bR^n)$ $\exists \e_0>0$ $\forall \e<\e_0$: $\sup_{x \in K}\abso{Q(S_\e\varphi, x)}<\frac{1}{2}$ and thus $S(S_\e\varphi, x)>0$. Furthermore, $\exists N \in \bN$ $\forall \varphi \in \cA_N(\bR^n)$ $\exists \e_1>0$ $\forall \e<\e_1$: $\sup_{x \in K}\abso{S(S_\e\varphi, x)}<\e^{-N}$. Then for $q \coleq \max(p,N)$, $\varphi \in \cA_q(\bR^n)$, $\e < \min(\e_0, \e_1)$ and $x \in K$ we obtain
\begin{align*}
\abso{R(S_\e\varphi, x)} &= \abso{\frac{1+Q(S_\e\varphi, x)}{S(S_\e\varphi, x)}} \ge \frac{1-\abso{Q(S_\e\varphi, x)}}{S(S_\e\varphi, x)} > \frac{\e^N}{2}.
\end{align*}

Conversely, given $R$ satisfying (ii) set $S(\varphi) \coleq 1/R(\varphi)$ where defined and $0$ elsewhere. Then $S \in \cE^e_M(\Omega)$ by definition and obviously $RS-1 \in \cN^e(\Omega)(n)$.
\end{proof}

The following proposition establishes a relation between invertibility and point values.

\begin{proposition} $\widetilde R \in \cG^e(\Omega)$ is invertible if and only if $\widetilde R(\widetilde X)$ is invertible in $\widetilde\bC(n)$ for each $\widetilde X \in \widetilde\Omega_c$.
\end{proposition}
\begin{proof}Necessity holds because point evaluation at a fixed generalized point evidently is a ring homomorphism from $\cG^e(\Omega)$ into $\widetilde\bC(n)$, thus $\widetilde R \widetilde S=1$ in $\cG^e(\Omega)$ implies $\widetilde R(\widetilde X)\widetilde S(\widetilde X)=1$ in $\widetilde\bC(n)$.
For sufficiency suppose that $\widetilde R$ is not invertible. Then by Theorem \ref{geinv} $\exists K \csub \Omega$ $\forall q \in \bN$ $\exists \varphi_q \in \cA_q(\bR^n)$ $\exists (\e_{q,k})_{k \in \bN} \searrow 0$ $\exists (x_{q,k})_{k \in \bN} \in K^\bN$ such that $\abso{R(S_{\e_{q,k}}\varphi_q, x_{q,k})} \le \e_{q,k}^q$.
Let $X: \cA_0(\bR^n) \rightarrow K$ and $(q_l)_{l \in \bN}$ be as obtained from Lemma \ref{constructpoint} with arbitrary $x_0 \in K$.
Then clearly $X$ is compactly supported and the class of $R(X)$ is not strictly nonzero and thus not invertible, because for arbitrary $q$ we can choose any $l$ such that $q_l\ge q$ and for large enough $k \in \bN$ we get
\begin{equation*}
\abso{R(X)(S_{\e_{q_l,k}}\varphi_{q_l})} = \abso{R(S_{\e_{q_l,k}}\varphi_{q_l}, x_{q_l, k})} \le \e_{q_l,k}^{q_l} \le \e^q_{q_l, k}.
\end{equation*}

\end{proof}

Proposition \ref{inverti} also follows directly from the following Lemma, whose validity is clear because for $\tilde r \in \widetilde \bC(n)$ and $\widetilde X \in \widetilde \Omega_c(n)$ we have $\tilde \rho(\tilde r)(\widetilde X) = \tilde r$.

\begin{lemma}$\tilde r \in \widetilde \bC(n)$ is invertible if and only if $\tilde \rho(\tilde r) \in \cG^e(\Omega)$ is.
\end{lemma}

\section{Point values in $\cG^d(\Omega)$}

While in $\cG^s$ and $\cG^e$ one can essentially leave away the $x$-slot in order to obtain generalized numbers we have to be more careful when introducing generalized numbers in the diffeomorphism invariant setting. First, smoothness of the involved objects is a crucial factor requiring considerable technical machinery (cf.~\cite{GKOS} Chapter 2). Second, there are two equivalent formalisms for describing the algebra $\cG^d$: one stems from the original construction by J.-F.~Colombeau \cite{ColElem}, the other is used by J.~Jel\'inek \cite{Jelinek} and is essential if one aims to construct a corresponding algebra intrinsically on a manifold. It is a sensible requirement that the translation mechanism between the C-formalism and the J-formalism (\cite{GKOS} Section 2.3.2) remains intact in order to translate results related to point values.

As we are dealing with differential algebras we can define generalized numbers as constant generalized functions, which means those functions $R$ satisfying $D_iR = 0$ $\forall i$. For connected $\Omega$ this is a natural definition of a space of numbers, generalized points simply are vectors of such numbers. Now as $D_i$ only acts on the $x$-slot one would be tempted to simply leave it away as we did in the $\cG^e$-setting with the hope to get simpler objects. We refrain from doing so, however, as viewing generalized numbers as a subspace of the generalized functions has two significant advantages: first, the existing technical background regarding smoothness which lies at the basis of $\cG^d$ can be used. Second, the translation mechanism given by the map $T$ works straightforward.

Instead of requiring $D_iR=0$ one can equivalently demand that the function does not depend on the second slot. We thus obtain the following definition.

\begin{definition} Let $V \subseteq \bR^p$ be open for some $p \in \bN$. Then generalized points of $V$ in the $\cG^d$-setting are defined by
\begin{align*}
\begin{split}
V_M(\Omega) &\coleq \{ X \in C^\infty(U(\Omega), V)\ |\ \forall K \csub \Omega\ \forall \alpha \in \bN_0^n\ \exists N\in \bN\\
&\qquad \forall \phi \in C^\infty_b(I \times \Omega, \cA_0(\bR^n)):\\
&\qquad \sup_{x \in K}\abso{\pd^\alpha X(S_\e\phi(\e,x),x)} = O(\e^{-N})\\
&\qquad \text{and }
\forall (\varphi, x),(\varphi ,y) \in U(\Omega):\ X(\varphi,x) = X(\varphi,y)\},
\end{split}\\
\begin{split}
V_N(\Omega) &\coleq \{ X \in V_M(\Omega)\ |\ \forall K \csub \Omega\ \forall m \in \bN\ \exists q \in \bN\\
&\qquad \forall \phi \in C^\infty_b(I \times \Omega,\cA_q(\bR^n)): \sup_{x \in K}\abso{X(S_\e\phi(\e,x),x)} = O(\e^m) \},
\end{split}\\
\widetilde V(\Omega) &\coleq V_M(\Omega) / V_N(\Omega).\\
\end{align*}
In order to obtain moderateness estimates of generalized point values one needs to introduce the concept of compactly supported generalized points, as is exemplified in the special algebra resp.~elementary full algebra by
\[ \abso{(u(x))_\e} = \abso{u_\e(x_\e)} \le \sup_{x \in K}\abso{u_\e(x)} \]
resp.~
\[ \abso{R(X)(S_\e\varphi)} = \abso{R(S_\e\varphi, X(S_\e\varphi))} \le \sup_{x \in K}\abso{R(S_\e\varphi,x)} \]
where $x_\e \in K$ for small $\e$ resp.~$X(S_\e\varphi) \in K$ for all $\varphi$ with sufficiently many vanishing moments and small $\e$. In order to find an analogous condition for $\cG^d$ one could start with a representative $X \in V_M(\Omega)$ of a generalized point satisfying $X(\varphi,x)\in L$ for all $(\varphi, x) \in U(\Omega)$ and some compact set $L \csub \Omega$. However, this condition is not preserved under change of representative: if one adds an element $Y$ of $V_N(\Omega)$ to $X$ one can only retain
\begin{gather*}
  \forall K \csub \Omega\ \exists q \in \bN\ \forall \phi\in C^\infty_b(I \times \Omega, \cA_q(\bR^n))\ \exists \e_0>0\\
  \forall \e<\e_0\ \forall x \in K: (X+Y)(S_\e\phi(\e,x),x) \in L'
\end{gather*}
where $L'$ is an arbitrarily small compact neighborhood of $L$.
The reason for this is that negligibility of $Y \in V_N(\Omega)$ gives uniformly small values of $Y(S_\e\phi(\e,x),x)$ (for $x \in K$ and $\e$ small) only if $\phi$ is an element of $C^\infty_b(I \times \Omega, \cA_q(\bR^n))$ for some certain $q$. This means that if $\phi$ has less than $q$ vanishing moments $Y(S_\e\phi(\e,x),x)$ may grow in any moderate way, leaving no hope of staying near $L$ or even in any compact subset of $V$, in general.

The easiest remedy to this problem is to simply define a generalized point $\widetilde X\in\widetilde V(\Omega)$ as being compactly supported if it has at least one representative $X$ whose image is contained in some compact set and only use such a suitable representative for the definition of point evaluation.

A different approach which is not pursued here but has to be mentioned is to use an equivalent description of $\cG^d(\Omega)$ where tests for moderateness and negligibility are performed using test objects having asymptotically vanishing moments. Such an algebra, called $\cG^2(\Omega)$, exists and is diffeomorphism invariant \cite{found}. It was demonstrated by J.~Jel\'inek \cite{Jelinek} that this algebra actually is the same as $\cG^d(\Omega)$. Using the moderateness and negligibility conditions of $\cG^2(\Omega)$ it would be possible to redefine the spaces used here in order to have a definition of compact support which is stable under change of representatives. In order to be consistent with our formalism of $\cG^d$, however, we chose not to take this route here, as it has no effect on the validity of the point value characterization theorem below and because there is no straightforward interface between $\cG^2(\Omega)$ and $\cG^d(\Omega)$.

\begin{definition}
A generalized point $\widetilde X \in \widetilde V(\Omega)$ is called compactly supported in $L \csub V$ if it has a representative $X \in V_M(\Omega)$ such that $\forall (\varphi, x) \in U(\Omega)$: $X(\varphi,x) \in L$. Denote by $\widetilde V_c(\Omega)$ the subset of all compactly supported generalized points of $\widetilde V(\Omega)$.
\end{definition}

As usual, elements of $V_M(\Omega)$ resp.~$V_N(\Omega)$ are called moderate resp.~negligible and we write $X \sim Y$ for $X-Y \in V_N(\Omega)$.
\end{definition}
Setting $V = \bC$ gives the space $\widetilde \bC(\Omega)$ of generalized complex numbers over $\Omega$.
As $X \in C^\infty(U(\Omega), V)$ is moderate resp.~negligible if and only if each component $\pr_i \circ X$ is, \cite{GKOS} Theorems 2.5.3 and 2.5.4 immediately give a characterization of moderateness resp.~negligibility of $X$ in terms of differentials of $X_\e \coleq X \circ S^{(\e)}$: $X \in C^\infty(U(\Omega), V)$ is moderate if and only if $\forall K \csub \Omega$ $\forall \alpha \in \bN_0^n$ $\forall k \in \bN_0$ $\exists N \in \bN$ $\forall B \subseteq \ccD(\bR^n)$ bounded it holds that
\begin{equation*}
\norm{\pd^\alpha\ud_1^k X_\e(\varphi,x)(\psi_1,\dotsc,\psi_k)} = O(\e^{-N}) \qquad (\e \rightarrow 0)
\end{equation*}
resp.~$X \in V_M(\Omega)$ is negligible if and only if $\forall K \csub \Omega$ $\forall m \in \bN$ $\exists q \in \bN$ $\forall B \subseteq \ccD(\bR^n)$ bounded it holds that
\begin{equation*}
\norm{X_\e(\varphi, x)} = O(\e^m) \qquad (\e \rightarrow 0),
\end{equation*}
where the estimate has to hold uniformly for $x \in K$, $\varphi \in B \cap \cA_0(\bR^n)$ resp.~$B \cap \cA_q(\bR^n)$, and $\psi_1,\dotsc,\psi_k \in B \cap \cA_{00}(\bR^n)$.

In the Colombeau-setting the point value is obtained as in $\cG^s$ and $\cG^e$ by inserting the (generalized) point into the $x$-slot. The corresponding formula for the J-setting is obtained by using the translation mechanism provided by the map $T^*$. We fix the following abbreviations for the natural definitions of point evaluation in the Jel\'inek- and the Colombeau-setting, noting that no confusion can arise from using the expression $R(X)$ in both cases.
\begin{enumerate}
  \item $R(X)(\varphi, x) \coleq R(T_{X(\varphi, x)-x}\varphi, X(\varphi, x))$ for $R \in C^\infty(\cA_0(\Omega) \times \Omega)$ and $X \in C^\infty(\cA_0(\Omega) \times \Omega, \Omega)$, and
  \item $R(X)(\varphi, x) \coleq R(\varphi, X(\varphi, x))$ for $R \in C^\infty(U(\Omega))$ and $X \in C^\infty(U(\Omega), \Omega)$.
\end{enumerate}

Because $R(X)$ is not defined on the whole of $\cA_0(\Omega) \times \Omega$ resp.~$U(\Omega)$, one has to implement a smooth cut-off procedure as in the following proposition. We will do so first in the Jel\'inek-setting because there the smoothness issues are more perspicuous -- the topology on $U(\Omega)$ is induced by the mapping $T$, so questions of smoothness on $U(\Omega)$ are most easily handled by transferring them to $\cA_0(\Omega) \times \Omega$.

\begin{proposition}\label{makeitsmooth}
Given $R \in C^\infty(\cA_0(\Omega) \times \Omega)$ and $X \in C^\infty(\cA_0(\Omega) \times \Omega, \Omega)$ satisfying
\begin{equation}\label{montag}
\exists L\csub \Omega\ \forall (\varphi, x) \in \cA_0(\Omega)\times \Omega: X(\varphi, x) \in L
\end{equation}
there exists a map $J_{R,X} \in C^\infty(\cA_0(\Omega) \times \Omega)$ such that for any $K\csub \Omega$ and any $B \subseteq \ccD(\bR^n)$ satisfying $\exists \beta>0$ $\forall \omega\in B$: $\supp \omega \subseteq \overline{B_\beta(0)}$ there is a relatively compact open neighborhood $U$ of $K$ in $\Omega$ and $\e_0>0$ such that for all $x \in U$, $\varphi \in B \cap \cA_0(\bR^n)$, and $\e < \e_0$ the expression $R(X)(T_xS_\e\varphi, x)$ is defined and 
\begin{equation*}
  J_{R,X}(T_xS_\e\varphi,x) = R(X)(T_xS_\e\varphi, x).
\end{equation*}
\end{proposition}

\begin{proof}
  Let $z \in \Omega$ remain fixed for the following construction. For some $\delta_z>0$ smaller than $\frac{1}{3} \dist(L, \pd\Omega)$ and $\frac{1}{2}\dist(z, \pd \Omega)$ we set $A_z \coleq B_{\delta_z}(z) \subseteq \Omega$ and $B_z \coleq B_{\delta_z}(A_z) = B_{2 \delta_z}(z)$. Both sets are relatively compact in $\Omega$. For all $x \in \overline{A_z}$ and $\varphi \in \cA_{0,\overline{B_z}}(\Omega)$ we consequently obtain
\begin{align*}
\supp T_{X(\varphi, x)-x} \varphi &= X(\varphi, x)-x + \supp \varphi\\
&\subseteq L-x + \overline{B_{2\delta_z}(z)} \subseteq L + \overline{B_{3\delta_z}(0)} \subseteq \Omega
\end{align*}
which means that $R(X)(\varphi,x) = R(T_{X(\varphi, x)-x}\varphi, X(\varphi, x))$ is defined on the set $\cA_{0,\overline{B_z}}(\Omega) \times A_z$. Furthermore $g_z \coleq R(X)|_{\cA_{0, \overline{B_z}}(\Omega) \times A_z} \in C^\infty(\cA_{0, \overline{B_z}}(\Omega) \times A_z)$: this follows easily by writing down all maps and spaces involved, after which $g_z$ is seen to be a composition of smooth functions. Set $D_z \coleq B_{\delta_z/2}(A_z)$ and choose a smooth function $\rho_z \in C^\infty(\Omega, \bR)$ with support in $\overline{B_z}$ and $\rho_z \equiv 1$ on $\overline{D_z}$. Fixing an arbitrary $\varphi_z \in \cA_{0,\overline{B_z}}(\Omega)$ define the projection
\[ \pi_z(\varphi) \coleq \varphi \cdot \rho_z + (1 - \int \varphi \cdot \rho_z) \cdot \varphi_z\qquad \forall \varphi \in \cA_0(\Omega), \]
then clearly $\pi_z \in C^\infty(\ccD(\bR^n), \ccD(\bR^n))$ and thus $\pi_z \in C^\infty(\cA_0(\Omega), \cA_{0, \overline{B_z}}(\Omega))$: the restriction to a set carrying the initial smooth structure with respect to the inclusion evidently is smooth, and as $\pi_z$ has values in $\cA_{0,\overline{B_z}}(\Omega)$ and this set also carries the initial smooth structure, $\pi_z$ is smooth into this set. For $\supp \varphi \subseteq \overline{D_z}$ we have $\pi_z(\varphi) = \varphi$.
There exists a smooth partition of unity $\{\chi_z\}_z$ subordinate to $\{A_z\}_z$, that is a collection of maps $\chi_z \in C^\infty(\Omega, [0,1])$ with $\supp \chi_z \subseteq A_z$ such that set of supports $\{\supp \chi_z\}_z$ is locally finite and $\sum \chi_z(x) =1$ $\forall x \in \Omega$.
Define a map $f_z$ on $\cA_0(\Omega) \times \Omega$ by
\begin{equation*}
f_z(\varphi,x) \coleq \left\{ \begin{split}
g_z(\pi_z(\varphi),x)\chi_z(x) & \qquad \text{if }x \in A_z\\
0 & \qquad \text{otherwise.}
\end{split}\right.
\end{equation*}
We have $f_z \in C^\infty(\cA_0(\Omega) \times \Omega)$: given an arbitrary smooth curve $c=(c_1,c_2) \in C^\infty(\bR, \cA_0(\Omega) \times \Omega)$, $c^*f_z$ is smooth, as any $t_0 \in \bR$ has a neighborhood whose image under $c_2$ lies either in $A_z$ or in the complement of $\supp \chi_z$, which are open sets covering $\Omega$. In the first case, \[ g_z(\pi_z(c_1(t)), c_2(t)) \chi_z(c_2(t)) = g_z(\pi_z(c_1(t)), \tilde c_2(t))\chi_z(\tilde c_2(t))\]
 in a neighborhood of $t_0$ on which $c_2$ is equal to some curve $\tilde c_2 \in C^\infty(\bR, A_z)$, thus one can employ smoothness of $(\varphi, x) \mapsto g_z(\pi_z(\varphi), x)\chi_z(x)$ on $\cA_0(\Omega) \times A_z$. In the second case the function is zero on an open neighborhood of $c_2(t_0)$, thus smooth trivially.
Now we can define $J_{R,X}: \cA_0(\Omega) \times \Omega \rightarrow \bC$ as $J_{R,X}(\varphi, x) \coleq \sum_z f_z(\varphi, x)$, which also is easily seen to be smooth as the sum is locally finite in $x$.
Now let $K$ and $B$ be given as stated in the proposition. $K$ has an open neighborhood $U$ which meets only finitely many supports of the $\chi_z$, which means that there are $z_1,\dotsc,z_m \in \Omega$ for some $m \in \bN$ such that $K \subseteq U \subseteq \bigcup_{i=1\dotsc m}\supp \chi_{z_i} \subseteq \bigcup_{i=1\dotsc m}A_{z_i}$, so on $\cA_0(\Omega) \times U$ $J_{R,X}$ is given by $\sum_{i=1\dotsc m}f_{z_i}$.
For $\e < \min_i \delta_{z_i}/(2\beta)$, $\varphi \in B \cap \cA_0(\bR^n)$ and $x \in A_{z_i}$, $\supp T_xS_\e\varphi \subseteq \overline{B_{\e\beta}(x)} \subseteq \overline{D_{z_i}}$ and thus $\pi_{z_i}(T_xS_\e\varphi) = T_xS_\e\varphi$; now $x \in \supp \chi_{z_i} \subseteq A_{z_i}$ implies $g_{z_i}(\pi_{z_i}(T_xS_\e\varphi), x) = R(X)(T_xS_\e\varphi,x)$ and thus for $x \in U$, $\varphi \in B \cap \cA_0(\bR^n)$, and $\e$ as above we finally obtain the conclusion
\begin{align*}
  J_{R,X}(T_xS_\e\varphi,x) & = \sum_{i=1\dotsc m}g_{z_i}(\pi_{z_i}(T_xS_\e\varphi),x)\chi_{z_i}(x)\\
&= R(X)(T_xS_\e\varphi,x) \sum_{i=1\dotsc m}\chi_{z_i}(x) = R(X)(T_xS_\e\varphi,x).\qedhere
\end{align*}
\end{proof}

\begin{corollary}\label{makeitsmoothC}
Given $R \in C^\infty(U(\Omega))$ and $X \in C^\infty(U(\Omega), \Omega)$ satisfying
\begin{equation}\label{montagabend}
\exists L \csub \Omega\ \forall (\varphi, x) \in U(\Omega):\ X(\varphi, x) \in L
\end{equation}
there exists $S_{R,X} \in C^\infty(U(\Omega))$ such that for any $K\csub \Omega$ and $B \subseteq \ccD(\bR^n)$ satisfying $\exists \beta>0$ $\forall \omega \in B$: $\supp \omega \subseteq \overline{B_\beta(0)}$ there is a relatively compact open neighborhood $U$ of $K$ in $\Omega$ and $\e_0>0$ such that for all $x \in U$, $\varphi \in B \cap \cA_0(\bR^n)$, and $\e < \e_0$, the expression $R(X)(S_\e\varphi, x)$ is defined and 
\begin{equation*}
S_{R,X}(S_\e\varphi,x) = R(X)(S_\e\varphi, x).
\end{equation*}
\end{corollary}
\begin{proof}
  We define $R^J \coleq (T^{-1})^*R \in C^\infty(\cA_0(\Omega) \times \Omega)$ and $X^J \coleq (T^{-1})^*X \in C^\infty(\cA_0(\Omega) \times \Omega, \Omega)$. Then $X^J$ satisfies \eqref{montag}, giving $J_{R^J, X^J} \in C^\infty(\cA_0(\Omega) \times \Omega)$. Now by Proposition \ref{makeitsmooth} there exists a relatively compact open neighborhood $U$ of $K$ in $\Omega$ and $\e_0>0$ such that $\forall x \in U$, $\varphi \in B \cap \cA_0(\bR^n)$, and $\e<\e_0$ we know that $R^J(X^J)(T_xS_\e\varphi, x)$ is defined and $J_{R^J, X^J}(T_xS_\e\varphi, x) = R^J(X^J)(T_xS_\e\varphi, x)$. Thus because $T^*(R^J(X^J)) = R(X)$ we obtain the result by setting $S_{R,X} \coleq T^*J_{R^J, X^J}$.
\end{proof}

The following proposition establishes that the construction of $S_{R,X}$ defines a unique element of $\widetilde \bC(\Omega)$ and enables us to use it for the definition of point values in $\cG^d(\Omega)$.

\begin{proposition}\label{srxprop}Given $R \in \cE^d_M(\Omega)$ and $X,Y \in \Omega_M(\Omega)$ satisfying \eqref{montagabend} $S_{R,X}$ is in $\bC_M(\Omega)$; if $R$ is negligible $S_{R,X}$ is, and $X \sim Y$ implies $S_{R,X} \sim S_{R,Y}$.
\end{proposition}

\begin{proof}
  Fix $K \csub \Omega$, $\alpha \in \bN_0^n$, and $k \in \bN_0$ for testing and let $B \subseteq \ccD(\bR^n)$ be bounded for testing in terms of differentials. Moderateness of $S_{R,X}$ is tested by estimating
\begin{equation*}
\abso{\pd^\alpha\ud_1^k(S_{R,X})_\e(\varphi,x)(\psi_1,\dotsc,\psi_k)}
\end{equation*}
where $x \in K$, $\varphi \in B \cap \cA_0(\bR^n)$, and $\psi_1,\dotsc,\psi_k \in B \cap \cA_{00}(\bR^n)$. Let $J \subseteq \bR$ be a bounded neighborhood of $0$. Then $B+J\psi_1+\dotsm+J\psi_k$ is bounded in $\ccD(\bR^n)$. Corollary \ref{makeitsmoothC} gives an open neighborhood $U$ of $K$ in $\Omega$ and $\e_0>0$ such that for $x \in U$, $\varphi \in B' \cap \cA_0(\bR^n)$, and $\e < \e_0$ the equation
\begin{equation*}
(S_{R,X})_\e (\varphi, x) = (R(X))_\e (\varphi, x)
\end{equation*}
holds. Given $\varphi,\psi_1,\dotsc,\psi_k$ as above we obtain for the $k$th differential
\begin{align*}
\ud_1^k(S_{R,X})_\e&(\varphi,x)(\psi_1,\dotsc,\psi_k) \\
&= \left.\frac{\pd}{\pd t_1}\right|_0\dotsm \left.\frac{\pd}{\pd t_k}\right|_0 (S_{R,X})_\e(\varphi + t_1\psi_1 + \dotsc + t_k\psi_k, x) \\
&= \left.\frac{\pd}{\pd t_1}\right|_0\dotsm \left.\frac{\pd}{\pd t_k}\right|_0 (R(X))_\e(\varphi + t_1\psi_1 + \dotsc + t_k\psi_k, x) \\
&= \ud_1^k ( R(X) )_\e(\varphi, x)(\psi_1,\dotsc,\psi_k).
\end{align*}

Note that this seemingly trivial equality and the following application of the chain rule rest on two hidden details. First, because in the first slot the mappings $S_{R,X}$ and $R(X)$ are defined on subsets of the affine subspace $\cA_0(\Omega)$, their differentials have to be calculated by considering the corresponding maps on the linear subspace $\cA_{00}(\Omega)$ which are obtained by pullback along an affine isomorphism $\cA_{00}(\Omega) \rightarrow \cA_0(\Omega)$. Second, these maps obtained actually have to be restricted to suitable subsets of $\cA_0(\bR^n) \times \Omega$ in order to give meaning to their differentials (cf.~\cite{GKOS} Section 2.3.3 for a detailed discussion).

As $(R(X))_\e(\varphi, x) = R_\e(\varphi, X_\e(\varphi,x))$, by the chain rule (\cite{GKOS} Appendix A) $\ud_1^k (R(X))_\e(\varphi, x)(\psi_1,\dotsc,\psi_k)$ consists of terms of the form 
\begin{equation*}
(\ud_2^l\ud_1^mR_\e(\varphi, X_\e(\varphi,x))(\psi_{i_1},\dotsc,\psi_{i_m}, (\ud_1^{a_1}X_\e)(\varphi, x)(\psi_{A_1}),\dotsc, (\ud_1^{a_l}X_\e(\varphi, x)(\psi_{A_l}))
\end{equation*}
where $m,l \in \bN_0$,
$i_1, \dotsc, i_m \in \{1\dotsc k\}$, $a_1,\dotsc,a_l \in \bN$, and $\psi_{A_1},\dotsc,\psi_{A_l}$ are appropriate tuples of elements from $\{\psi_1,\dotsc,\psi_k\}$. Consequently, the expression $\pd^\alpha\ud_1^k(R(X))_\e(\varphi,x)(\psi_1,\dotsc,\psi_k)$ consists of terms of the form
\begin{align*}
(\ud_2^l\ud_1^m\pd^\gamma R_\e(\varphi, X_\e&(\varphi,x)))(\psi_{i_1},\dotsc,\psi_{i_m},\\
&\pd^{\beta_1}\ud_1^{a_1}X_\e(\varphi, x)(\psi_{A_1}),\dotsc,\pd^{\beta_l}\ud_1^{a_l}X_\e(\varphi, x)(\psi_{A_l}))
\end{align*}
where $\gamma, \beta_1,\dotsc,\beta_l$ are some multi-indices. The norm of the last expression can be estimated by
\begin{multline*}
  \norm{(\ud_2^l\ud_1^m \pd^\gamma R_\e(\varphi, X_\e(\varphi,x))(\psi_{i_1},\dotsc,\psi_{i_m})}\cdot\\
\cdot \norm{(\pd^{\beta_1}\ud_1^{a_1}X_\e(\varphi, x)(\psi_{A_1})} \dotsm \norm{(\pd^{\beta_l}\ud_1^{a_l}X_\e(\varphi, x)(\psi_{A_l})}
\end{multline*}
whence the first two claims of the proposition follow immediately from moderateness resp.~negligibility of $R$ and moderateness of the compactly supported $X$.

For the last claim, fix $K \csub \Omega$ and $m \in \bN$ for testing and let $B \subseteq \ccD(\bR^n)$ be bounded. Let $Y$ take values in $L \csub \Omega$. We need to estimate $\abso{ (S_{R,X} - S_{R,Y})_\e(\varphi, x)}$ for $x \in K$ and $\varphi \in B \cap \cA_0(\bR^n)$. By Corollary \ref{makeitsmoothC} there exists an open neighborhood $U$ of $K$ in $\Omega$ such that for $x \in U$, $\varphi \in B \cap \cA_0(\bR^n)$, and small $\e$ we have both $(S_{R,X})_\e(\varphi, x) = (R(X))_\e(\varphi,x)$ and $(S_{R,Y})_\e(\varphi, x) = (R(Y))_\e(\varphi, x)$, so we have to estimate $\abso{(R(X) - R(Y))_\e(\varphi, x)}$.
Setting $F(t) \coleq R_\e(\varphi, (Y + t(X-Y))_\e(\varphi, x))$ the last expression can be written as $\abso{F(1)-F(0)}$.
As $X \sim Y$ there exists $q \in \bN$ such that for $x \in K$, $\varphi \in B \cap \cA_q(\bR^n)$, and small $\e$ we have $\abso{(X-Y)_\e(\varphi, x)} < \e$, so $F(t)$ is defined and smooth on $[0,1]$ and we can write

\begin{equation*}
\abso{F(1)-F(0)} = \abso{\int_0^1 F'(t)\ud t} = \abso{\int_0^1\ud_2R_\e(\varphi, g(t)) \cdot (X- Y)_\e(\varphi,x) \ud t}
\end{equation*}
whence the claim follows directly from moderateness of $R$ and negligibility of $X-Y$.
\end{proof}

\begin{definition}
  For $\widetilde R \in \cG^d(\Omega)$ and $X \in \widetilde \Omega_c(\Omega)$ we define the generalized point value of $\widetilde R$ at $\widetilde X$ as $\widetilde R (\widetilde X) \coleq [S_{R,X}]$ where $R$ is any representative of $\widetilde R$ and $X$ is a representative of $\widetilde X$ satisfying \eqref{montagabend}.
\end{definition}

\begin{lemma}\label{accpoint}Let $K$ be a compact set. Given for each $q \in \bN$ a sequence $(x_{q,k})_{k \in \bN}$ in $K$ it holds that
\begin{multline*}
\exists x_0 \in K\ \forall \delta>0\ \forall q_0\in\bN\ \exists q=q(\delta,q_0)\ge q_0\\
\forall k_0 \in \bN\ \exists k=k(\delta, q_0, k_0)\ge k_0:\ x_{q,k} \in B_{\delta}(x_0).
\end{multline*}
This means that $x_0$ is an accumulation point of infinitely many of the sequences $(x_{q,k})_k$.
\begin{proof}
Assuming the converse we would have $\forall x_0 \in K$ $\exists \delta=\delta(x_0) > 0$ $\exists q_0 = q_0(x_0) \in \bN$ $\forall q \ge q_0$ $\exists k_0 = k_0(x_0, q)$ $\forall k \ge k_0$: $x_{q,k} \not\in B_{\delta(x_0)}(x_0)$. As $K \subseteq \bigcup_{x \in K}B_{\delta(x)}(x)$ we can choose $x_1,\dotsc,x_m$ $(m \in \bN)$ such that $K$ is contained in $\bigcup_{i=1,\dotsc,m}B_{\delta(x_i)}(x_i)$. Then for $q \ge \max_i q_0(x_i)$ and $k \ge \max_i k_0(x_i, q)$ we obtain the contradiction $x_{q,k} \not\in \bigcup_{i=1,\dotsc,m}B_{\delta(x_i)}(x_i) \supseteq K$.
\end{proof}
\end{lemma}

After these preparations we are finally able to establish the point value characterization theorem for $\cG^d(\Omega)$.

\begin{theorem}\label{gdpvchar}$\widetilde R \in \cG^d(\Omega)$ is $0$ if and only if $\widetilde R(\widetilde X)=0$ in $\widetilde \bC(\Omega)$ for all $\widetilde X \in \widetilde \Omega_c(\Omega)$.
\end{theorem}

\begin{proof}
  Let $R$ be a representative of $\widetilde R$.
  We have already shown in Proposition \ref{srxprop} that $R \in \cN^d(\Omega)$ implies $R(X) \in \bC_\cN(\Omega)$ for all $X \in \Omega_M(\Omega)$. For the converse we assume $R \not\in \cN^d(\Omega)$ and construct a generalized point $X$ such that $R(X) \not\in \bC_\cN(\Omega)$. By this assumption there exists $K \csub \Omega$ and $m \in \bN$ such that for all $q \in \bN$ there is some $\phi_q \in C^\infty_b(I \times \Omega, \cA_q(\bR^n))$ such that $\forall k \in \bN$ $\exists \e_{q,k}<\frac{1}{k}$ $\exists x_{q,k} \in K$ such that with $\varphi_{q,k} \coleq S_{\e_{q,k}}\phi_q(\e_{q,k}, x_{q,k})$ we have
\[ \abso{R(\varphi_{q,k}, x_{q,k})} \geq \e_{q,k}^m. \]
For the negligibility test of $R(X)$ to fail it suffices to construct $X$ such that for each of infinitely many $q$ the equation $X(\varphi_{q,k}, x_{q,k}) = x_{q,k}$ holds for infinitely many $k$.
Choose positive real numbers $\delta$ and $\eta_1$ both smaller than $\dist(x_0, \pd\Omega)$. Lemma \ref{accpoint} gives
\begin{equation}\label{afds}
\begin{aligned}
\exists x_0\in K\ \forall q_0\in\bN\ &\exists q=q(\delta,q_0)\ge q_0\ \forall k_0 \in \bN\\
& \exists k=k(\delta, q_0, k_0)\ge k_0: x_{q,k} \in B_{\delta}(x_0).
\end{aligned}
\end{equation}

Furthermore, for all $q \in \bN$ there exists an index $k_1(q) \in \bN$ such that $\supp S_{\e_{q,k}} \phi_q (\e_{q,k}, x_{q,k}) \subseteq B_{\eta_1}(0)$ for all $k \ge k_1(q)$. Combining this with \eqref{afds}, there exists a strictly increasing sequence $(q_l)_{l \in \bN}$ and for each $l\in \bN$ a sequence $(k_{l,r})_{r \in \bN}$ with $k_{l,r} \ge k_1(q_l)$ and $x_{q_l, k_{l,r}} \in B_\delta(x_0)$ for all $r \in \bN$. Choose $\eta_2>0$ arbitrary and set $U \coleq \{ \varphi \in \ccD(\bR^n)\ |\ \norm{\varphi}_\infty < \eta_2 \}$.

Let $(c_n)_{n \in \bN}$ be a sequence in $\bN$ in which each natural number appears infinitely often.
Set $\varphi_1 \coleq \varphi_{q_{c_1}, k_{c_1,1}}$ and $x_1 \coleq x_{q_{c_1}, k_{c_1,1}}$. Inductively, given $\varphi_n$ choose $r$ large enough such that $\norm{\varphi_{q_{c_{n+1}}, k_{c_{n+1},r}}}_\infty > \norm{\varphi_n}_\infty + 2\eta_2$ and set $\varphi_{n+1} \coleq \varphi_{q_{c_{n+1}}, k_{c_{n+1},r}}$ and $x_{n+1} \coleq x_{q_{c_{n+1}}, k_{c_{n+1},r}}$.

The sequences $(\varphi_n)_{n \in \bN}$ and $(x_n)_{n \in \bN}$ then have the following properties:
\begin{enumerate}
\item $x_n \in B_{\delta}(x_0)\ \forall n \in \bN$.
\item For each of infinitely many $q \in \bN$ there are infinitely many $k \in \bN$ such that $\varphi_{q,k}$ resp.~$x_{q,k}$ appears in the sequence $(\varphi_n)_n$ resp.~$(x_n)_n$.
\item $\supp \varphi_n \subseteq B_{\eta_1}(0)$ for all $n \in \bN$.
\item Because $\norm{\varphi_n}_\infty = \norm{T_{-x_n}\varphi_n}_\infty$ all sets $\overline{U} + T_{-x_n}\varphi_n$ for $n \in \bN$ are pairwise disjoint.
\end{enumerate}

Choose $\eta_3$ such that $0<\eta_3 < \eta_2$. Set $U' \coleq \{ \varphi \in \ccD(\bR^n)\ |\ \norm{\varphi}_\infty <\eta_3 \}$, $\bE \coleq \ccD_{\overline{B_{\eta_1}(0)}}(\bR^n)$ and $U'_1 \coleq U' \cap \bE$.
Construct a smooth bump function $\chi_1 \in C^\infty(\bE, \bR)$ with $\supp \chi_1 \subseteq \overline{U_1'}$ and $\chi_1(0)=1$ as follows:

Let $g \in C^\infty(\bR, \bR)$ be nonnegative such that $g(x)=1$ for $x \leq 0$ and $g(x)=0$ for $x \geq 1$. As $\bE$ is a nuclear locally convex space,
there exist a convex, circled $0$-neighborhood $V \subseteq U_1'$ and a positive semi-definite sesquilinear form $\sigma$ on $\bE$ such that $p: x \mapsto \sqrt{\sigma(x,x)}$ is the gauge function of $V$ and a continuous seminorm on $\bE$
(\cite{Schaefer} Chapter III 7.3).
From the Cauchy-Schwartz inequality we infer
$\abso{\sigma(x,y)} \leq p(x)p(y)$, which means that $\sigma$ is bounded and thus smooth. Consequently the associated hermitian form $h: x \mapsto \sigma(x,x)$ also is smooth.
The differentials of $h$ are given by
\begin{align*}
\ud h (x)(v) &= 2 \Re \sigma(x, v), \\
\ud^2 h (x)(v, w) &= 2 \Re \sigma(v,w),\textrm{ and} \\
\ud^3 h &= 0
\end{align*}
where $\Re$ denotes the real part. Now $\chi_1 \coleq g \circ h$ is in $C^\infty(\bE, \bR)$ with $\chi_1(0) = 1$ and $\supp \chi_1 \subseteq \overline{V} \subseteq \overline{U_1'} \subseteq U \cap \bE$ because $g(h(x)) >0$ implies $h(x)<1$ and thus $x \in V$.

Then by \cite{KM} Lemma 16.6 and an obvious adaptation of the proof of \cite{KM} Proposition 16.7 there exists a function $\chi \in C^\infty(\ccD(\bR^n), \bR)$ such that $\chi|_{\bE} = \chi_1$, $\chi(0)=1$ and $\supp \chi \subseteq \overline{U}$.

Set $\chi_m(\varphi) \coleq \chi(\varphi - T_{-{x_m}}\varphi_m)$  for $\varphi \in \ccD(\bR^n)$. We define a map $Y: \ccD(\bR^n) \times \bR^n \rightarrow \Omega$ by
\[ Y(\varphi,x) \coleq \sum_{m \in \bN} \left(x_0 + \chi_m(T_{-x}\varphi)(x_m-x_0)\right) \in B_\delta(x_0). \]
Because the supports of $\chi_m$ are disjoint $Y$ has at most one summand near any given $\varphi$; it clearly is smooth and as $\cA_0(\Omega) \times \Omega$ carries the initial smooth structure with respect to the inclusion its restriction to $\cA_0(\Omega) \times \Omega$ also is smooth. Our prospective generalized point is defined as
\[ X \coleq T^* ( Y|_{\cA_0(\Omega) \times \Omega}) \in C^\infty(U(\Omega), \Omega), \]
and satisfies $X(\varphi_n, x_n) = x_n$. $X$ is compactly supported in $\overline{B_\delta(x_0)}$.  
In order to show moderateness of $X$ we test in terms of differentials. Fix $K \csub \Omega$, $\alpha \in \bN_0^n$, $k \in \bN_0$, and $B \subseteq \ccD(\bR^n)$ bounded for testing. We then need to estimate the expression
\begin{equation*}
\pd^\alpha\ud_1^kX_\e(\varphi,x)(\psi_1,\dotsc,\psi_k)
\end{equation*}
where $x \in K$, $\varphi \in B \cap \cA_0(\Omega)$, and $\psi_1,\dotsc,\psi_k \in B \cap \cA_{00}(\Omega)$. We first look at the function whose derivatives we need:
\[ X_\e(\varphi, x) = Y(T_xS_\e\varphi, x) = \sum_m(x_0 + \chi_m(S_\e\varphi)(x_m-x_0)). \]
As we see from the right hand side this expression does not depend on $x$ so we only need to consider the case $\alpha=0$. If the $k$th differential at $\varphi$ in directions $\psi_1,\dotsc,\psi_k$ is nonzero it is given by only one term of the right hand side, so for each $\varphi$ there exists an index $m_0\in \bN$ such that 

\begin{multline*}
\ud_1^kX_\e(\varphi,x)(\psi_1,\dotsc,\psi_k) \\
= \ud^k \bigl( \varphi \mapsto (x_0 + \chi_{m_0}(S_\e\varphi)(x_{m_0}-x_0)) \bigr) (\varphi)(\psi_1, \dotsc, \psi_k) \\
= \ud^k \bigl( \varphi \mapsto (x_0 + \chi(S_\e\varphi - T_{-x_{m_0}}\varphi_{m_0})(x_{m_0}-x_0)) \bigr) (\varphi)(\psi_1, \dotsc, \psi_k)
\end{multline*}

In order to use that $\chi|_{\bE} = \chi_1$ we need that the support of the argument of $\chi$ in the previous expression is contained in $\overline{B_{\eta_1}(0)}$. By construction this is the case for all $\varphi_n$ and if $\e$ is small enough it is also satisfied for $S_\e\varphi$ for all $\varphi \in B$ uniformly. As $\chi_1 = g \circ h$ we need to obtain the differentials
\begin{equation}\label{dienstag}
  \ud^k \bigl( \varphi \mapsto g(h(S_\e\varphi - T_{-x_{m_0}}\varphi_{m_0})) \bigr) (\varphi, x)(\psi_1, \dotsc, \psi_k).
\end{equation}

Abbreviate $f(\varphi) \coleq S_\e\varphi - T_{-x_{m_0}}\varphi_{m_0}$. We can assume that $h(f(\varphi)) < 1$ holds, as otherwise expression \eqref{dienstag} vanishes. By the chain rule we see that the $k$th differential is given by the product of derivatives of $g$ (which are globally bounded) and terms of the form $\ud^k(h \circ f)(\varphi)(\psi_1,\dotsc,\psi_k)$ for some $k \in \bN$ which again by the chain rule are given by terms of the form
\begin{equation}\label{dienstaggamma}
(\ud^k h)(f(\varphi))(\ud^{l_1}f(\varphi)(\psi_{A_1}),\dotsc, \ud^{l_k}f(\varphi)(\psi_{A_k}))
\end{equation}
for some $l_1,\dotsc,l_k \in \bN$ and appropriate subsets $\psi_{A_1},\dotsc,\psi_{A_k} \subseteq \{ \psi_1,\dotsc,\psi_k \}$. Here only $k=0,1,2$ are relevant as higher derivatives of $h$ vanish. We obtain from \eqref{dienstaggamma} the three terms
\begin{align*}
h(f(\varphi)) &= \sigma(f(\varphi), f(\varphi))\\
\ud h(f(\varphi))(\ud f(\varphi)(\psi_1)) &= 2 \Re \sigma (f(\varphi), \ud f(\varphi)(\psi_1))\\
\ud^2 h(f(\varphi))(\ud f(\varphi)(\psi_1), \ud f(\varphi)(\psi_2)) &= 2 \Re \sigma (\ud f(\varphi)(\psi_1), \ud f(\varphi)(\psi_2))
\end{align*}
The function $f$ is differentiated at most once because its higher order derivatives vanish. Noting that $\ud f(\varphi)(\psi) = S_\e\psi$ we estimate these terms by the Cauchy-Schwartz inequality. We obtain products of
$\sqrt{h(f(\varphi))}$
(which has been assumed to be smaller than $1$) and
$\sqrt{h(S_\e\psi)} = p(S_\e\psi)$
(where $\psi$ is $\psi_1$ or $\psi_2$). Being a continuous seminorm, $p$ is majorized by finitely many of the usual seminorms $q_\alpha$ of $\bE$ given by $q_\alpha(\varphi) = \sup_{x \in \bR^n}\abso{\pd^\alpha\varphi(x)}$ for all $\alpha \in \bN_0^n$. We thus end up with the expression
\begin{align*}
  q_\alpha(S_\e\psi) &= \sup_{x \in \bR^n} \abso{\pd^\alpha (S_\e\psi)(x)} = \sup_{x \in \bR^n}\abso{\pd^\alpha(\e^{-n}\psi(x/\e))}\\
  &= \sup_{x \in \bR^n}\abso{\e^{-n-\abso{\alpha}}(\pd^\alpha\psi)(x/\e)} =  \e^{-n-\abso{\alpha}} \norm{\psi}_\infty
\end{align*}
and as $\psi$ is from the bounded set $B$ we get the desired growth estimates independently of $m_0$ and conclude that $X$ is moderate. By construction $R(X)$ is not negligible and the point value characterization theorem is established.
\end{proof}

\begin{acknowledgement}
The author thanks Michael Kunzinger for helpful discussions. This research has been supported by START-project Y237 and project P20525 of the Austrian Science Fund and the Doctoral College 'Differential Geometry and Lie Groups' of the University of Vienna. 
\end{acknowledgement}

\providecommand{\WileyBibTextsc}{}
\let\textsc\WileyBibTextsc
\providecommand{\othercit}{}
\providecommand{\jr}[1]{#1}
\providecommand{\etal}{~et~al.}

\end{document}